\documentclass [12pt]{article}
\usepackage {multicol}

\usepackage{amssymb,amsmath}
\usepackage{hyperref}
\usepackage{graphicx}
\usepackage{color}
\usepackage{caption}

\begin{document}

\title{\Large\bf Deleting Items and Disturbing Mesh Theorems for Riemann Definite Integral and Their Applications}

\author{\normalsize LIU Jingwei \footnote{\scriptsize{Corresponding Author: Email: liujingwei03@tsinghua.org.cn (J.W Liu); liuy0605@126.com (Y. Liu) } } $^{,1}$,\  LIU Yi $^{2}$ \\
{\scriptsize{(1. School of Mathematics and System Sciences, Beihang University, Beijing,100083P.R.China)}}\\
{\scriptsize{(2. School of Mathematics and Statistics, Beijing Institute of Technology,Beijing 100081 P.R.China)}}
}

\date{\normalsize February 12,2017}
\maketitle

\textbf{Abstract:} Based on the definition of Riemann definite integral,deleting items and disturbing mesh theorems on Riemann sums are given. After deleting some items or disturbing the mesh of partition, the limit of Riemann sums still converges to Riemann definite
integral under specific conditions. These theorems can deal with a class of complicate limitation of sum and product of series of a function, and demonstrate that the geometric intuition of Riemann definite integral is more profound than ordinary thinking of area of curved trapezoid.


\section{Introduction}

It is well known that mathematical analysis is an advanced mathematical branch founded on limitation. All definitions of limitation, continuity,
derivative, definite integral, improper integrals, infinite series, limitation of function of several variables, partial derivatives, multiple
integral, line integral, surface integral, parametric variable integral, \textit{etc,} are based on limitation definition. Definite integral, which is strictly defined on limitation by Bernhard Riemann (1826-1864) , adopts classical
four steps (partitioning interval, choosing distinguished point, making
Riemann sum, determining limit of Riemann sum) to approximate the area of
curved trapezoid, this rigorous framework of definition make it one of the
core research of differential and integral calculus [1-10]. Meanwhile, in
recent more than 150 years, the geometric meaning of Riemann definite
integral with area of curved trapezoid is commonly accepted in mathematical
field.

However, there is a theoretical question that only under the framework of
Riemann definite integral that the Riemann sum can approach the area of
curved trapezoid? This Riemann definition is the only way to calculate the
area of curved trapezoid? Two categories of theorems are developed in this
paper to show that under some proper conditions, deleting items and
disturbing mesh of Riemann sum, we still can obtain the limit of Riemann sum
with Riemann definite integral, which demonstrate that the geometric meaning
of Riemann definite integral is beyond area of curved trapezoid.

\section{Definition of definite integral and Deleting item \& Disturbing mesh Theorems}

Riemann definite integral can be described as follows [1-10]:

\textbf{Definition 1.} Suppose that $f(x)$ is a real function defined on a closed interval $[a,b]$, partition$\{x_0 ,x_1 , \cdots ,x_n \}$, such that $a = x_0 < x_1 < x_2 < \cdots < x_n = b$, divides $[a,b]$ into $n$ sub-intervals $[x_{k - 1} , x_k ](1 \le k \le n)$. Let
$\Delta x_k = x_k - x_{k-1}, (1 \le k \le n)$ and $\lambda =
\mathop {\max} \limits_{1 \le k \le n} \Delta x_k $, $\lambda $is called
mesh of the largest width of $n $sub-intervals. On each subinterval
$[x_{k - 1} ,x_k ]$, a distinguished point $\xi _k \in [x_{k - 1} ,x_k ]$ is
adopted and the Riemann sums is defined as $\sigma = \sum\limits_{k = 1}^n
{f(\xi _k )\Delta x_k}  $.

If for any partition and tagged point$\xi _k $, Riemann sums $\sigma $ converges to a common limit$ I$as$\lambda \to 0$, that is

\[
\mathop {\lim} \limits_{\lambda \to 0} \sigma = \mathop {\lim
}\limits_{\lambda \to 0} \sum\limits_{k = 1}^n {f(\xi _k )} \Delta x_k = I
\]

We call$ I$is the Riemann definite integral of$f(x)$on $[a,b]$, and
denote $I = \int_a^b {f(x)dx} $. Also, $f(x)$is called Riemann integrable
on$ [a,b]$.

In short, the $\varepsilon-\delta$ definition for Riemann integral is as follows:

$\forall \varepsilon > 0, \exists \delta > 0,$for any partition of $[a,b]$, if $\lambda = \mathop {\max}\limits_{1 \le k \le n} \Delta x_k < \delta $, for any$\xi _k \in [x_{k -1} ,x_k ]$, such that $\vert \sigma - I\vert < \varepsilon $. Then, $\mathop {\lim} \limits_{\lambda \to 0} \sum\limits_{k = 1}^n {f(\xi _k )} \Delta x_k = I = \int_a^b {f(x)dx} $.

For convenience of discussion, we give the following definition.

\textbf{Definition 2.} Suppose that $f(x)$is a real function defined on a closed interval $[a,b]$, partition$\{x_0 ,x_1 , \cdots ,x_n \}$ such that $a = x_0 < x_1 < x_2 < \cdots < x_n = b $divides$ [a,b]$ into $n$sub-intervals $[x_{k - 1} x_k ](1 \le k \le n)$. Let $\Delta x_k = x_k -
x_{k-1}, (1 \le k \le n)$ and $\lambda = \mathop {\max} \limits_{1\le k \le n} \Delta x_k $.On each $[x_{k - 1} ,x_k ]$, a distinguished
point$\xi _k \in [x_{k - 1} ,x_k ]$ is adopted. Denote $J = \{1,2, \cdots,n\}$. For a fixed natural number$K \in N^ + (0 < K < n)$, denote$J_K =
\{i_1 , \cdots ,i_K \} \subset J$. And,

\begin{equation}
\label{eq1}
\bar {\sigma} _K = \sum\limits_{k \in J / J_K}  {f(\xi _k )\Delta x_k}
\end{equation}

The problem of calculating $\mathop {\lim} \limits_{\lambda \to 0} \bar{\sigma}_K $ is called the deleting item problem in Riemann definite
integral.

\textbf{Theorem 1. (Deleting item Theorem)} Suppose that$ f(x)$is Riemann integrable on$[a,b]$.Partition$\{x_0 ,x_1 , \cdots ,x_n \}$ such
that $a = x_0 < x_1 < x_2 < \cdots < x_n = b$ divides$ [a,b]$ into $n $sub-intervals $[x_{k - 1} x_k ](1 \le k \le n)$. Let $\Delta x_k
= x_k - x_{k - 1} (1 \le k \le n)$ and$\lambda = \mathop {\max} \limits_{1\le k \le n} \Delta x_k $.For any distinguished point $\xi _k $,$x_{k - 1}\le \xi _k \le x_k $, and a fixed natural number$K \in N^ + (0 < K < n)$,the deleting item Riemann sum is $\bar {\sigma} _K = \sum\limits_{k \in J /J_K}  {f(\xi _k )\Delta x_k}  $.Then,

\begin{equation}
\label{eq2}
\mathop {\lim} \limits_{\lambda \to 0} \sum\limits_{k \in J / J_K}  {f(\xi_k )\Delta x_k}  = I = \int_a^b {f(x)dx} .
\end{equation}

\textbf{Proof:} Since $f(x)$is Riemann integrable on $[a,b]$, $f(x)$ is bounded on$ [a,b]$, that is$ \exists M > 0$,such that $\vert
f(x)\vert \le M$. And,

$\forall \varepsilon > 0, \exists \delta = \frac{\varepsilon} {2MK} > 0,$for any partition of $[a,b]$, if$\lambda = \mathop {\max} \limits_{1 \le k \le n} \Delta x_k < \delta $,for any distinguished point$\xi _k \in [x_{k - 1} ,x_k ]$,we have

\[
\vert \sum\limits_{k \in J} {f(\xi _k )\Delta x_k}  - I\vert <\frac{\varepsilon}{2}.
\]

Again,

\[
\vert \sum\limits_{k \in J_K}  {f(\xi _k )\Delta x_k}  \vert \le MK\lambda <\frac{\varepsilon}{2}.
\]

Then,
\[
\begin{array}{ll}
 \vert \sum\limits_{k \in J / J_K}  {f(\xi _k )\Delta x_k}  - I\vert &= \vert\sum\limits_{k \in J} {f(\xi _k )\Delta x_k}  - I - \sum\limits_{k \in J_K}{f(\xi _k )\Delta x_k}  \vert \\
      & \le \vert \sum\limits_{k \in J} {f(\xi _k )\Delta x_k}  - I\vert + \vert \sum\limits_{k \in J_K}  {f(\xi _k )\Delta x_k} \vert
      <\displaystyle \frac{\varepsilon}{2}+\frac{\varepsilon}{2}=\varepsilon \\
 \end{array}
\]

That is $\mathop {\lim} \limits_{\lambda \to 0} \sum\limits_{k \in J / J_K}
{f(\xi _k )\Delta x_k}  = I$. Hence, it ends the proof. \hfill $\Box$

For convenience of description, if a special case is considered that
interval $[a,b]$ is partitioned into $n$ equal width sub-interval, and $\xi _k $is adopted right endpoint of sub-interval. We have theorem 2.

\textbf{Theorem 2. (Deleting Item Theorem)} Suppose that$ f(x)$ is Riemann integrable on $[a,b]$, $K$ is a given natural number. Then

\begin{equation}
\label{eq3}
\mathop {\lim} \limits_{n \to + \infty}  \frac{b - a}{n}\sum\limits_{k \in J/ J_K}  {f(a + \frac{k}{n}(b - a))} = I.
\end{equation}

Specially, if$J_K = \{1,2, \cdots ,K\}$,then

\begin{equation}
\label{eq4}
\mathop {\lim} \limits_{n \to + \infty}  \frac{b - a}{n}\sum\limits_{k = K +
1}^n {f(a + \frac{k}{n}(b - a))} = I.
\end{equation}

If $J_K = \{n - K + 1,n - K + 2, \cdots ,n\}$, then

\begin{equation}
\label{eq5}
\mathop {\lim} \limits_{n \to + \infty}  \frac{b - a}{n}\sum\limits_{k =1}^{n - K} {f(a + \frac{k}{n}(b - a))} = I.
\end{equation}

Instinctively, the geometric meaning of Riemann definite integral is using the sum of areas of rectangles to approach the area of curved trapezoid.
However, theorem 1 and theorem 2 show that in the framework of limitation, abandoning finite items of Riemann sum (infinitesimal), the residue sum
still converges to area of curved trapezoid. Though K would be very large, it is a fixed number, abandoning the area of $K$ rectangle is just the
infinitesimal in process of $\lambda \to 0$.

Deleting item theorem of definite integral demonstrate that the nature of Riemann definite integral is far more beyond the simply comprehension of
area of curved trapezoid. Next, we will give another generalization theorem of Riemann sums. For convenience, we only give the special case of partition with $n$ equal sub-interval length of $[a,b]$, and distinguished point adopts right endpoint of sub-interval , which is $\xi _k = x_k $.

\textbf{Definition 3}.Suppose that partition $\{x_0 ,x_1 , \cdots ,x_n \}$, such that $a = x_0 < x_1 < x_2 < \cdots < x_n = b$, divides
$[a,b]$into $n$ sub-intervals. Let $\tilde {\Delta} x_k $denote a disturbed value of $\Delta x_k = x_k - x_{k - 1} $, namely, $\tilde {\Delta} x_k $ has a small amount of change of $\Delta x_k $. Specially, in the equal sub-interval width case of partition,$\lambda = \frac{b - a}{n}$, let
$d_{n,k} ( \ge 0)(1 \le k \le n)$ denote a non-negative disturbance , we call

\begin{equation}
\label{eq6}
\tilde {\Delta} x_k = \frac{b - a}{n + d_{n,k}}
\end{equation}

as the disturbance of $\Delta x_k $, and sequence $d_{n,k}$ is the disturbance sequence of sequence $\Delta x_k $.

The disturbance sequence in definition 3 perturbs the sub-interval width of partition of $[a,b]$ in Riemann definite integral, hence it disturbs mesh of partition. And, $\tilde {\Delta} x_k - \Delta x_k $ is expected to be infinitely small.

\textbf{Theorem 3. (Disturbing Mesh Theorem)} Suppose that $f(x)$ is non-negative and Riemann integrable on$[a,b]$. Given an equal sub-interval
width partition of $[a,b]$, let $d_{n,k} (1 \le k \le n)$ be the non-negative disturbance sequence of $\Delta x_k $. If

\begin{equation}
\label{eq7}
\mathop {\lim} \limits_{n \to + \infty}  \frac{d_{n,k}} {n} = 0 ,
1 \le k \le n.
\end{equation}

Then,

\begin{equation}
\label{eq8}
\mathop {\lim} \limits_{n \to + \infty}  \sum\limits_{k = 1}^n {f(a +\frac{k}{n}(b - a))\frac{b - a}{n + d_{n,k}} } = I = \int_a^b {f(x)dx} .
\end{equation}

\textbf{Proof:} Denote $d_{\min}  = \mathop {\min} \limits_{1 \le k \le n}\{d_{n,k} \}$, $d_{\max}  = \mathop {\max} \limits_{1 \le k \le n} \{d_{n,k}\}$.

Then, $\mathop {\lim} \limits_{n \to + \infty}\displaystyle  \frac{n}{n + d_{\min}}=1$, $\mathop {\lim} \limits_{n \to + \infty} \displaystyle \frac{n}{n + d_{\max}}=1$.

Since,
\[
\sum\limits_{k = 1}^n {\frac{f(a + \frac{k}{n}(b - a))(b - a)}{n + d_{\max} }} \le \sum\limits_{k = 1}^n {\frac{f(a + \frac{k}{n}(b - a))(b - a)}{n +d_k} } \le \sum\limits_{k = 1}^n {\frac{f(a + \frac{k}{n}(b - a))(b - a)}{n+ d_{\min} } }
\]

\[
\mathop {\lim} \limits_{n \to + \infty}  \sum\limits_{k = 1}^n {\frac{f(a +\frac{k}{n}(b - a))(b - a)}{n + d_{\max} } } = \mathop {\lim} \limits_{n \to+ \infty}  \frac{n}{n + d_{\max} } \sum\limits_{k = 1}^n {f(a +\frac{k}{n}(b - a))\frac{(b - a)}{n}}=I
\]

\[
\mathop {\lim} \limits_{n \to + \infty}  \sum\limits_{k = 1}^n {\frac{f(a +\frac{k}{n}(b - a))(b - a)}{n + d_{\min} } } = \mathop {\lim} \limits_{n \to+ \infty}  \frac{n}{n + d_{\min} } \sum\limits_{k = 1}^n {f(a +\frac{k}{n}(b - a))\frac{(b - a)}{n}}= I
\]

According to squeeze theorem, we obtain
\[
\mathop {\lim} \limits_{n \to + \infty}  \sum\limits_{k = 1}^n {\frac{f(a +\frac{k}{n}(b - a))(b - a)}{n + d_{n,k}} } = I= \int_a^b {f(x)dx} .
\]
\hfill $\Box$

\textbf{Remark 1:} A concise condition for$\mathop {\lim} \limits_{n \to +\infty}\displaystyle  \frac{d_{n,k}} {n} = 0 $,$1 \le k \le n$, is as
follows:

$\exists \quad p < 1,\quad p \in R$, and $d_{n,k} = O(n^p)$.    \hfill (9)

\textbf{Corollary 1. (Disturbing Mesh Theorem)} Suppose that $f(x)$ is non-negative and Riemann integrable on$[a,b]$. Given an equal sub-interval
width partition of $[a,b]$, let $d_{n,k} (1 \le k \le n)$ be the non-negative disturbance sequence of $\Delta x_k $. If

$\exists \quad p < 1, p \in R,$and $d_{n,k} = O(n^p)$. Then,

\[
\mathop {\lim} \limits_{n \to + \infty}  \sum\limits_{k = 1}^n {f(a +\frac{k}{n}(b - a))\frac{b - a}{n + d_{n,k}} } = I = \int_a^b {f(x)dx} .
\]

Combining Theorem 2 and Theorem 3, we have Theorem 4:

\textbf{Theorem 4. (Deleting item-Disturbing mesh Theorem)} Suppose that$f(x)$is non-negative and Riemann integrable on$[a,b]$. Given a partition
of $[a,b]$, $d_{n,k} (1 \le k \le n)$ is the non-negative disturbance sequence of $\Delta x_k $.$K$ is a given natural number. If

\[
\mathop {\lim} \limits_{n \to + \infty}  \frac{d_{n,k}} {n} = 0 , 1 \le k \le n.
\]

Then,
\[
\mathop {\lim} \limits_{n \to + \infty}  \sum\limits_{k \in J / J_K}{\frac{f(a + \frac{k}{n}(b - a))(b - a)}{n + d_{n,k}} } = I = \int_a^b
{f(x)dx} .    \tag{10}
\]

Specially, If$J_K = \{1,2, \cdots ,K\}$, then

\[
\mathop {\lim} \limits_{n \to + \infty}  \sum\limits_{k = K + 1}^n{\frac{f(a + \frac{k}{n}(b - a))(b - a)}{n + d_{n,k}} } = I = \int_a^b
{f(x)dx} .    \tag{11}
\]

If $J_K = \{n - K + 1,n - K + 2, \cdots ,n\}$, then
\[
\mathop {\lim} \limits_{n \to + \infty}  \sum\limits_{k = 1}^{n - K}{\frac{f(a + \frac{k}{n}(b - a))(b - a)}{n + d_{n,k}} } = I = \int_a^b
{f(x)dx} .    \tag{12}
\]

\textbf{Remark 2.} The condition of $f(x)$is non-negative and integrable on $[a,b]$ in Theorem 3 and Theorem 4 can be relaxed to $f(x)$ is non-negative or non-positive and integrable on $[a,b]$.And, if $f(x)$can be decomposed into finite non-negative or non-positive sub-intervals, and the conditions of Theorem 3 and Theorem 4 hold on each sub-interval, so do the corresponding conclusions. However, for the whole interval $[a,b]$, the method of squeeze theorem may not be applied, we omit the discussion.

\textbf{Remark 3.} Theorem 3 and Theorem 4 can be also described in style of Theorem 1 in general form of Riemann definite integral definition, the
depiction is not concise, we omit it. Theorem 3 show that, in the common sense of definition of Riemann integral ,not only the traditional four-step Riemann sums can approach definite integral, but also the corresponding sums that the width sequence is disturbed by a high order infinitely small
disturbance of $\lambda $, can reach the limit of definite integral. Obviously, the disturbance sequence can include a large amount of functions
and real sequences. From the conclusion of theorem 4, many of Riemann sums with deleting item and disturbing mesh can approach the same definite
integral value. Hence, using the definition of Riemann definite integral, a lot of complicated Riemann sums limitation problems can be solved.

\textbf{Remark}\textbf{ 4.} In the above discussion, we assume that $K$ is a fixed natural number, however, if we examine the proof of theorem 1, if $K$ is related to $\lambda $ (also $n)$, in some conditions, for example, $K(\lambda )$ satisfying

$K(\lambda ) > 0,$and $K(\lambda ) \cdot \lambda \to 0$,$(\lambda \to 0)$.                      \hfill (13)

Theorem 1 still holds. Thus, in Theorem 2 and Theorem 4, while $K$ is related to $n$ and satisfies the following condition:

$K(n) > 0$,and $\displaystyle \frac{K(n)}{n} \to 0, n \to + \infty. $     \hfill (14)

Theorem 2 and Theorem 4 can also hold. A concise case for this condition is:

$K(n)>0$ and  $\exists 0 < q <1,   q \in R$, such that $K(n) = O(n^q)$    \hfill (15)

Therefore, for Theorem 2 and 4, we have the following extensions:

\textbf{Theorem 5. (Deleting item Theorem)} Suppose that $f(x)$ is integrable on $[a,b]$. And $K(n)$ is a natural number, satisfying $K(n) >
0$,
\[
\frac{K(n)}{n} \to 0,n \to + \infty .
\]
Then,
\[
\mathop {\lim} \limits_{n \to + \infty}  \frac{b - a}{n}\sum\limits_{k \in J/ J_{K(n)}}  {f(a + \frac{k}{n}(b - a))} = I = \int_a^b {f(x)dx} . \tag{16}
\]

Specially, if $J_{K(n)} = \{1,2, \cdots ,K(n)\}$, then

\[
\mathop {\lim} \limits_{n \to + \infty}  \frac{b - a}{n}\sum\limits_{k =K(n) + 1}^n {f(a + \frac{k}{n}(b - a))} = I = \int_a^b {f(x)dx} . \tag{17}
\]

If $J_{K(n)} = \{n - K(n) + 1,n - K(n) + 2, \cdots ,n\}$, then

\[
\mathop {\lim} \limits_{n \to + \infty}  \frac{b - a}{n}\sum\limits_{k =1}^{n - K(n)} {f(a + \frac{k}{n}(b - a))} = I = \int_a^b {f(x)dx} . \tag{18}
\]

\textbf{ Theorem 6. (Deleting item --Disturbing mesh Theorem)} Suppose that$f(x)$is non-negative and integrable on $[a,b]$. Given an equal
sub-interval width partition of $[a,b]$, suppose that $d_{n,k} (1 \le k\le n)$ is the non-negative disturbance sequence of sub-interval
width$\Delta x_k $, for a natural number $K(n) > 0$, if the following
conditions are satisfied

(1) $\mathop {\lim} \limits_{n \to + \infty}  \frac{d_{n,k}} {n} = 0 $,$1 \le k \le n$.

(2) $\frac{K(n)}{n} \to 0$,$n \to + \infty $

Then,

\[
\mathop {\lim} \limits_{n \to + \infty}  \sum\limits_{k \in J / J_{K(n)}} {\frac{f(a + \frac{k}{n}(b - a))(b - a)}{n + d_{n,k}} } = I = \int_a^b
{f(x)dx} .   \tag{19}
\]

Specially, if $J_{K(n)} = \{1,2, \cdots ,K(n)\}$, then

\[
\mathop {\lim} \limits_{n \to + \infty}  \sum\limits_{k = K(n) + 1}^n{\frac{f(a + \frac{k}{n}(b - a))(b - a)}{n + d_{n,k}} } = I = \int_a^b
{f(x)dx} .  \tag{20}
\]

If $J_{K(n)} = \{n - K(n) + 1,n - K(n) + 2, \cdots ,n\}$, then

\[
\mathop {\lim} \limits_{n \to + \infty}  \sum\limits_{k = 1}^{n - K(n)}{\frac{f(a + \frac{k}{n}(b - a))(b - a)}{n + d_{n,k}} } = I = \int_a^b
{f(x)dx} .    \tag{21}
\]

According to the concise conditions of formulae (9)(15), Theorem 5 and Theorem 6 have the following corollaries respectively.

\textbf{Corollary 2. (Deleting item Theorem)} Suppose that$f(x)$is integrable on $[a,b]$. And $K(n)$ is a natural number, satisfying$K(n) >
0$and $\exists      0< q < 1,   q \in R$, such that$K(n) =O(n^q)$. Then,

$\mathop {\lim} \limits_{n \to + \infty}  \frac{b - a}{n}\sum\limits_{k \in J/ J_{K(n)}}  {f(a + \frac{k}{n}(b - a))} = I = \int_a^b {f(x)dx}
.$

Specially, if $J_{K(n)} = \{1,2, \cdots ,K(n)\}$, then
\[
\mathop {\lim} \limits_{n \to + \infty}  \frac{b - a}{n}\sum\limits_{k =K(n) + 1}^n {f(a + \frac{k}{n}(b - a))} = I = \int_a^b {f(x)dx} .
\]

If $J_{K(n)} = \{n - K(n) + 1,n - K(n) + 2, \cdots ,n\}$, then
\[
\mathop {\lim} \limits_{n \to + \infty}  \frac{b - a}{n}\sum\limits_{k =1}^{n - K(n)} {f(a + \frac{k}{n}(b - a))} = I = \int_a^b {f(x)dx} .
\]

\textbf{Corollary 3. (Deleting item-Disturbing mesh Theorem)} Suppose that $f(x)$is non-negative and integrable on $[a,b]$. Given an equal
sub-interval width partition of $[a,b]$, suppose that $d_{n,k} (1 \le k\le n)$ is the non-negative disturbance sequence of sub-interval
width$\Delta x_k $, for a natural number $K(n) > 0$, if the following conditions are satisfied

(1) $\exists p < 1$,$p \in R$,and $d_{n,k} = O(n^p)$.

(2) $\exists      0 <q < 1,   q \in R$, such that$K(n) = O(n^q)$

Then,

\[
\mathop {\lim} \limits_{n \to + \infty}  \sum\limits_{k \in J / J_{K(n)}}{\frac{f(a + \frac{k}{n}(b - a))(b - a)}{n + d_{n,k}} } = I = \int_a^b
{f(x)dx} .
\]

Specially, if $J_{K(n)} = \{1,2, \cdots ,K(n)\}$,then

\[
\mathop {\lim} \limits_{n \to + \infty}  \sum\limits_{k = K(n) + 1}^n{\frac{f(a + \frac{k}{n}(b - a))(b - a)}{n + d_{n,k}} } = I = \int_a^b
{f(x)dx} .
\]

If $J_{K(n)} = \{n - K(n) + 1,n - K(n) + 2, \cdots ,n\}$, then

\[
\mathop {\lim} \limits_{n \to + \infty}  \sum\limits_{k = 1}^{n - K(n)}{\frac{f(a + \frac{k}{n}(b - a))(b - a)}{n + d_{n,k}} } = I = \int_a^b
{f(x)dx} .
\]

\textbf{Remark 5.} Obviously, the thought in this paper can be parallelly extended to Riemann definite integral with multiple variables, line integral and surface integral [11,12], there will exist the corresponding deleting item theorem and disturbing mesh theorems. We omit the concrete
descriptions. Meanwhile, In theorem 1-6, $\xi _k $ can be adopted as the left endpoint of $k$-th sub-interval, middle point, and any special point.
We also omit the verbose discussions.

\section{Application}

Two classical examples [8,13,14] are involved in discussion to show the applications of deleting item and disturbing mesh theorems.

\textbf{Example 1.}$^{[8]}$ $\mathop {\lim} \limits_{n \to \infty}  \sum\limits_{k = 1}^n {\displaystyle \frac{\sin (\frac{k}{n}\pi )}{n}} = \int_0^1
{\sin (x\pi )} dx = \frac{1}{\pi} \int_0^\pi {\sin (x)} dx = \frac{2}{\pi} $

\textbf{Example 1'.}$^{[13]}$ $\mathop {\lim} \limits_{n\to + \infty}  \sum\limits_{k = 1}^n {\displaystyle \frac{\sin \frac{k\pi} {n}}{n +
\frac{1}{k}} = \frac{1}{\pi} \int_0^\pi {\sin xdx} = \frac{2}{\pi} } $

\textbf{Example 2.} $^{[14]}$ Suppose that $f(x)$ is positive and integrable on $[0,1]$. Then,
\[
\mathop {\lim} \limits_{n \to \infty}  \sqrt[n]{f(\frac{1}{n})f(\frac{2}{n})\cdots f(\frac{n}{n})} =\displaystyle e^{ \mathop {\lim} \limits_{n \to \infty} \frac{1}{n}\sum\limits_{k = 1}^n {\ln f(\frac{k}{n})}}  =\displaystyle e^{\int_0^1 {\ln f(x)dx}}
\]

Based on the theory proposed in this paper, there will be following conclusions.

\subsection{Examples of deleting item theorem}

\ \ \ \ \ \textbf{Example 3.} $K$ is a fixed natural number, $\mathop {\lim} \limits_{n \to \infty}  \sum\limits_{k = K + 1}^n {\displaystyle \frac{\sin(\frac{k}{n}\pi )}{n}} = \int_0^1 {\sin (x\pi )} dx$

\textbf{Example 4.} $\mathop {\lim} \limits_{n \to \infty } \sum\limits_{k = [n^{\frac{8}{9}}] + 1}^n {\displaystyle \frac{\sin (\frac{k}{n}\pi
)}{n}} = \int_0^1 {\sin (x\pi )} dx$

\textbf{Example 5.} $\mathop {\lim} \limits_{n \to \infty } \sum\limits_{k = 1}^{n - [n^{\frac{8}{9}}]} {\displaystyle \frac{\sin (\frac{k}{n}\pi
)}{n}} = \int_0^1 {\sin (x\pi )} dx$

\textbf{Example 6.} Suppose that $f(x)$ is positive and integrable on $[0,1]$. Then,
\[
\mathop {\lim} \limits_{n \to \infty}  \sqrt[n]{f(\frac{K + 1}{n})f(\frac{K + 2}{n}) \cdots f(\frac{n}{n})}  = e^{\mathop {\lim
}\limits_{n \to \infty}  \frac{1}{n}\sum\limits_{k = K + 1}^n {\ln f(\frac{k}{n})}}  = e^{\int_0^1 {\ln f(x)dx}}
\]

\textbf{Example 7.} Suppose that $f(x)$ is positive and integrable on $[0,1]$. Then,
\[
\mathop {\lim} \limits_{n \to \infty}  \sqrt[n]{f(\frac{[n^{\frac{8}{9}}] + 1}{n})f(\frac{[n^{\frac{8}{9}}] + 2}{n}) \cdots f(\frac{n}{n})} = e^{\mathop {\lim} \limits_{n \to \infty}  \frac{1}{n}\sum\limits_{k = [n^{\frac{8}{9}}]+ 1}^n {\ln f(\frac{k}{n})}}  = e^{\int_0^1 {\ln f(x)dx}}
\]

\textbf{Example 8.} Suppose that $f(x)$ is positive and integrable on $[0,1]$. Then,
\[
\mathop {\lim} \limits_{n \to \infty}  \sqrt[n]{f(\frac{1}{n})f(\frac{2}{n})\cdots f(\frac{n - [n^{\frac{8}{9}}]}{n})} = e^{\mathop {\lim} \limits_{n\to \infty}  \frac{1}{n}\sum\limits_{k = 1}^{n - [n^{\frac{8}{9}}]} {\ln f(\frac{k}{n})}}  = e^{\int_0^1 {\ln f(x)dx}}
\]

\subsection{Examples of disturbing mesh theorem}

\ \ \ \ \ \textbf{Example 9.} $\mathop {\lim} \limits_{n \to \infty } \sum\limits_{k = 1}^n {\displaystyle \frac{\sin (\frac{k}{n}\pi )}{n + \frac{k}{n}}} = \int_0^1 {\sin (x\pi )} dx$

This example is also be proved by squeeze theorem, as $\displaystyle \sum\limits_{k = 1}^n {\frac{\sin (\frac{k}{n}\pi )}{n + 1}} \le
\sum\limits_{k = 1}^n {\frac{\sin (\frac{k}{n}\pi )}{n + \frac{k}{n}}} \le \sum\limits_{i = 1}^n {\frac{\sin (\frac{k}{n}\pi )}{n + \frac{1}{n}}}
$. According to Theorem 3, we can easily reach the conclusion.

\textbf{Example 10.} $\mathop {\lim} \limits_{n \to \infty } \sum\limits_{k = 1}^n {\displaystyle \frac{\sin (\frac{k}{n}\pi )}{n +
k^{[\frac{8}{9}]}}} = \int_0^1 {\sin (x\pi )} dx$

\textbf{Example 11.} Suppose that $f(x)$ is positive and integrable on $[0,1]$.Then,
\[
\mathop {\lim} \limits_{n \to \infty}  \sqrt[{n +[n^{\frac{8}{9}}]}]{f(\frac{1}{n})f(\frac{2}{n}) \cdots f(\frac{n}{n})}
=e^{\mathop {\lim} \limits_{n \to \infty}  \sum\limits_{k = 1}^n {\ln f(\frac{k}{n})} \frac{1}{n + [n^{\frac{8}{9}}]}}
= e^{\int_0^1 {\ln f(x)dx}}
\]

\textbf{Example 12.} Suppose that $f(x)$ is positive and integrable on $[0,1]$, and $\ln f(x)$ is non-negative or non-positive on
$[0,1]$.Then,
\[
\mathop {\lim} \limits_{n \to \infty}  \prod\limits_{k = 1}^n {[f(\frac{k}{n})]^{\frac{1}{n + \frac{k}{n}}}}  = e^{\mathop
{\lim} \limits_{n \to \infty}  \sum\limits_{k = 1}^n {\ln f(\frac{k}{n})} \frac{1}{n + \frac{k}{n}}} = e^{\int_0^1 {\ln f(x)dx}}
\]

In examples 10 -12, $\frac{1}{n + [k^{\frac{8}{9}}]}$, $\frac{1}{n+[n^{\frac{8}{9}}]}$ and $\frac{1}{n + \frac{k}{n}}$ are the little disturbances
of width of sub-intervals (also the mesh). Example 1' is also the case of disturbance.

\subsection{Examples of deleting item - disturbing mesh theorem}

\ \ \ \ \ \textbf{Example 13.} $\mathop {\lim} \limits_{n \to \infty}  \sum\limits_{k = [n^{\frac{8}{9}}] + 1}^n {\frac{\sin (\frac{k}{n}\pi )}{n
+ \frac{k}{n}}} = \int_0^1 {\sin (x\pi )} dx$

\textbf{Example 14.} $\mathop {\lim} \limits_{n \to \infty}  \sum\limits_{k = 1}^{n - [n^{\frac{8}{9}}]} {\frac{\sin (\frac{k}{n}\pi
)}{n + \frac{k}{n}}} = \int_0^1 {\sin (x\pi )} dx$

\textbf{Example 15.} $\mathop {\lim} \limits_{n \to \infty}  \sum\limits_{k = [n^{\frac{8}{9}}] + 1}^n {\frac{\sin (\frac{k}{n}\pi )}{n
+ \sqrt k} } = \int_0^1 {\sin (x\pi )} dx$

\textbf{Example 16.} $\mathop {\lim} \limits_{n \to \infty}  \sum\limits_{k = 1}^{n - [n^{\frac{8}{9}}]} {\frac{\sin (\frac{k}{n}\pi
)}{n + \sqrt k} } = \int_0^1 {\sin (x\pi )} dx$

\textbf{Example 17.} Suppose that $f(x)$ is positive and integrable on $[0,1]$. Then,
\[
\mathop {\lim} \limits_{n \to \infty}  \sqrt[{n + \sqrt n
}]{f(\frac{[n^{\frac{8}{9}}] + 1}{n})f(\frac{[n^{\frac{8}{9}}] + 2}{n})
\cdots f(\frac{n}{n})} = e^{\mathop {\lim} \limits_{n \to \infty
} \sum\limits_{k = [n^{\frac{8}{9}}] + 1}^n {\ln f(\frac{k}{n})\frac{1}{n +
\sqrt n} }}  = e^{\int_0^1 {\ln f(x)dx}}
\]

\textbf{Example 18.} Suppose that $f(x)$ is positive and integrable on $[0,1]$. Then,

$\mathop {\lim} \limits_{n \to \infty}  \sqrt[{n + \sqrt n
}]{f(\frac{1}{n})f(\frac{2}{n}) \cdots f(\frac{n - [n^{\frac{8}{9}}]}{n})} =
e^{\mathop {\lim} \limits_{n \to \infty}  \sum\limits_{k = 1}^{n -
[n^{\frac{8}{9}}]} {\ln f(\frac{k}{n})\frac{1}{n + \sqrt n} }}  =
e^{\int_0^1 {\ln f(x)dx}} _{ \newline
}$

\textbf{Example 19.} Suppose that $f(x)$ is positive and integrable on $[0,1]$, and $\ln f(x)$ is non-negative or non-positive on
$[0,1]$. Then,

\[
\mathop {\lim} \limits_{n \to \infty}  \prod\limits_{k = [n^{\frac{8}{9}}] +1}^n {[f(\frac{k}{n})]^{\frac{1}{n + \frac{k}{n}}}}
= e^{\mathop {\lim} \limits_{n \to \infty}  \sum\limits_{k =
[n^{\frac{8}{9}}] + 1}^n {\ln f(\frac{k}{n})} \frac{1}{n + \frac{k}{n}}} =
e^{\int_0^1 {\ln f(x)dx}}
\]

\textbf{Example 20.} Suppose that $f(x)$ is positive and integrable on $[0,1]$, and $\ln f(x)$ is non-negative or non-positive on
$[0,1]$. Then,

\[
\mathop {\lim} \limits_{n \to \infty}  \prod\limits_{k = 1}^{n -[n^{\frac{8}{9}}]} {[f(\frac{k}{n})]^{\frac{1}{n + \frac{k}{n}}}}  = e^{\mathop {\lim} \limits_{n \to \infty}  \sum\limits_{k = 1}^{n -[n^{\frac{8}{9}}]} {\ln f(\frac{k}{n})} \frac{1}{n + \frac{k}{n}}} = e^{\int_0^1 {\ln f(x)dx}}
\]

\section{Conclusion}

The deleting item and disturbing mesh theorems show that Riemann definite integral is not only the area of curved trapezoid, but also the limit of
magnitude sums in limitation theory. By arbitrary partitioning, choosing distinguished points, making up Riemann sums, only if the Riemann sums
approaches a limit as $\lambda \to 0$, it is defined as the Riemann definite integral. When some of the items are abandoned from Riemann sums (finite or infinite items) or the mesh has a little disturbance, only if the performances do not affect the trend of limit, this limit is still equal to the area of curved trapezoid. This mathematical meaning is different from the finite case of sums and the simple comprehension of Riemann definite integral's area of curved trapezoid background. Particularly, in real-world application, given a partition of $n$ size, abandoning finite items definitely affects the accuracy of calculation. However, in the sense of limitation, the theorems in this paper do not affect the area of curved trapezoid, since they are dependent on $\lambda \to 0$, also$n \to + \infty $. The results also show the gap of mathematical theory and
application, and furthermore, mathematical theory is far more complicated than simply geometric intuition. Even for the simple measurement problem of curved trapezoid, the conclusion can be reached that area of curved trapezoid can be calculated by Riemann definite integral, inversely, Riemann
definite integral can comprise magnitude sums besides Riemann sums.

Meanwhile, the research of this paper provides solutions of the limitation problems of a class of sum and product of function value satisfying deleting item and disturbing mesh theorems. This work is the generation and supplement of Riemann definite integral in high level sense of limitation and can also be extended to the Riemann definite integral of multiple variables, line integral and surface integral ,which are in the same framework of partitioning interval, choosing distinguished point, making up Riemann sums, and calculating limit. As there is relationship between Riemann definite integral and Lebesgue integral, so to integral and measure, the research in this paper indicates a more deep thought of mathematical nature of measure.

\section*{Acknowledgments}
This work is partially supported by Major Program of the National Natural Science Foundation of China (No.61327807).

\end{document}